\documentclass[reqno,12pt]{amsart}

\usepackage{vmargin}
\setpapersize{A4}

\usepackage{amsmath}

\usepackage{amsfonts}

\usepackage{amssymb}

\usepackage{eufrak}

\usepackage{amscd}

\usepackage{amsthm}

\usepackage{epsfig}

\usepackage{amstext}

\usepackage[all]{xy}
\CompileMatrices 

   \newcommand{\Hom}{\operatorname{Hom}}

\newcommand{\Ad}{\operatorname{Ad}}

\newcommand{\id}{\operatorname{id}}
\newcommand{\Aut}{\operatorname{Aut}}
\newcommand{\eq}[1]{\begin{equation}#1\end{equation}}
\newcommand{\eQ}[1]{\begin{equation*}#1\end{equation*}}
\newcommand{\spl}[1]{\begin{split}#1\end{split}}
 \newcommand{\Ext}{\operatorname{Ext}}
\newcommand{\diag}{\operatorname{diag}}

\newcommand{\cone}{\operatorname{cone}}

   \theoremstyle{plain}
   \newtheorem{thm}{Theorem}[section]
   
   \newtheorem{lemma}[thm]{Lemma}
   \newtheorem{cor}[thm]{Corollary}
   \theoremstyle{definition}

   \theoremstyle{remark}
   
   \newtheorem{remark}[thm]{Remark}

\title[The Connes-Higson construction is an isomorphism]{The Connes-Higson
construction is an isomorphism}
\author{Vladimir Manuilov and Klaus Thomsen}

\begin{document}

\begin{abstract}
Let $A$ be a separable $C^*$-algebra and $B$ a stable $C^*$-algebra
containing a strictly positive element. We show that the group $\Ext(SA,B)$
of unitary equivalence classes of extensions of $SA$ by $B$, modulo the
extensions which are asymptotically split, coincides with the group of
homotopy classes of such extensions. This is done by proving that the
Connes-Higson construction gives rise to an isomorphism between
$\Ext(SA,B)$ and the $E$-theory group $E(A,B)$ of homotopy classes of
asymptotic homomorphisms from $S^2A$ to $B$.

\end{abstract}

\maketitle

\section{Introduction}

The fundamental homotopy functors on the category of separable
$C^*$-algebras are all based on extensions --- either \it a priori \rm or
\it a posteriori \rm . So also the $E$-theory of Connes and Higson;
in the words of the founders: 'La $E$-theorie est ainsi le quotient par
homotopie de la th\'eorie des extensions', cf. \cite{CH}. The connection
between the asymptotic homomorphisms which feature explicitly in the
definition of $E$-theory, and $C^*$-extensions, appears as a fundamental
construction which associates an asymptotic homomorphism $SA \to B$ to a
given extension of $A$ by $B$. While it is easy to see that the homotopy
class of the asymptotic homomorphism only depends on the homotopy class
of the extension it is not so easy to decide if the converse is also true;
if the extensions must be homotopic when the asymptotic homomorphisms
which they give rise to via the Connes-Higson construction are. A part of
the main result in the present paper asserts that this is the case when
$A$ is a suspension and $B$ is stable. Rather unexpectedly it turned out
that the methods we developed for this were also able to characterize
$E$-theory as the quotient of all extensions of $SA$ by $B$ by an algebraic
relation which is very similar to the algebraic relation which has been
considered on the set of extensions since the way-breaking work of Brown,
Douglas and Fillmore, \cite{BDF}. Recall that in the BDF-approach two
$C^*$-extensions are identified when they become unitarily equivalent
after addition by extensions which are split, meaning that the quotient
map admits a $*$-homomorphism as a right-inverse. In the algebraic
relation, on the set of all $C^*$-extensions of $SA$ by $B$, which we
will show gives rise to $E$-theory, two extensions are identified when
they become unitarily equivalent after addition by extensions which are
asymptotically split, where we call an extension
\eQ{\begin{xymatrix}{
0 \ar[r]  &  B \ar[r]  & E  \ar[r]^-p  &  A \ar[r]  & 0 }
\end{xymatrix}}
asymptotically split when there is an asymptotic homomorphism
$\pi = (\pi_t)_{t \in [1,\infty)} : A \to E$ such that $p \circ \pi_t =
\id_{A}$ for all $t$. We emphasize that with this relation all extensions
of $SA$ by $B$ admit an inverse. In contrast, Kirchberg has shown,
\cite{Ki}, that the unitary equivalence classes of extensions of $SA$
by $\mathcal K$, modulo the split extensions, do not form a group when
$A$ is the reduced group $C^*$-algebra of a discrete non-amenable subgroup
of a connected Lie-group. Since our results show that the algebraic
relation we have just described is the same as homotopy, our main result
can also be considered as a result on homotopy invariance and it is
therefore noteworthy that the proof is self-contained, and in particular
does not depend on the homotopy invariance results of Kasparov.

Since there is also an equivariant version of $E$-theory, \cite{GHT},
which is being used in connection with the Baum-Connes conjecture, we
formulate and prove our results in the equivariant case. With the present
technology this does not require much additional work, but since some of
the material which we shall build on does not explicitly consider the
equivariant setting, notably \cite{DL} and \cite{H-LT}, there are a few
places where we leave the reader to check that the results from these
sources can be adapted to the equivariant case.

\section{An alternative to the BDF extension group}

Let $G$ be a locally compact, $\sigma$-compact group, and let $A$ and
$B$ be separable $G$-algebras, i.e. separable $C^*$-algebras with a
pointwise norm-continuous action of $G$ by automorphisms. Assume also
that $B$ is weakly stable, i.e. that $B$ is equivariantly isomorphic to
$B \otimes \mathcal K$ where $\mathcal K$ denotes the compact operators
of $l_2$ with the trivial $G$-action. Let $M(B)$ denote the multiplier
algebra of $B$, $Q(B) = M(B)/B$ the corresponding corona algebra and
$q_B : M(B) \to Q(B)$ the quotient map. Then $G$ acts by automorphisms
on both $M(B)$ and $Q(B)$\footnote{These actions are not pointwise
normcontinuous in general.}. It follows from \cite{Th1} that we can
identify the set of equivariant $*$-homomorphisms, $\Hom_G(A,Q(B))$,
from $A$ to $Q(B)$ with the set of $G$-extensions of $A$ by $B$. Two
$G$-extensions $\varphi,  \psi : A \to Q(B)$ are \it unitarily equivalent
\rm when there is a unitary $w \in M(B)$ such that $q_B(w) \in Q(B)$
is $G$-invariant and $ \Ad q_B(w) \circ \varphi = \psi$. Since $B$ is
weakly stable the set of unitary equivalence classes of extensions of
$A$ by $B$ form a semi-group; the addition is obtained by choosing two
$G$-invariant isometries $V_1,V_2 \in M(B)$ such that $V_1V_1^* + V_2V_2^*
= 1$ and setting $\varphi \oplus \psi = q_B(V_1)\varphi(\cdot) q_B(V_1)^*
+ q_B(V_2)\psi(\cdot) q_B(V_2)^*$. A $G$-extension $\varphi : A  \to Q(B)$
will be called \it asymptotically split \rm when there is an asymptotic
homomorphism $\pi = \{\pi_t\}_{t \in [1,\infty)} : A \to M(B)$ such that
$q_B \circ \pi_t = \varphi$ for all $t$.
All asymptotic homomorphisms we consider in this paper will be assumed
to be equivariant in the sense that $\lim_{t \to \infty} g \cdot
\pi_t(a) - \pi_t(g \cdot a) = 0$ for all $a \in A$ and $g \in G$.
As in \cite{MT2} we say that a $G$-extension $\varphi : A \to Q(B)$
is \it semi-invertible \rm when there is a $G$-extension $\psi \in
\Hom_G(A,Q(B))$ such that $\varphi \oplus \psi : A \to Q(B)$ is
asymptotically split. Two semi-invertible
extensions, $\varphi, \psi$, are called \it stably unitary equivalent
\rm when they become unitarily equivalent after addition by
asymptotically split extensions, i.e.
when there is an asymptotically split extension $\lambda $ such that
$\varphi \oplus \lambda $ is unitarily equivalent to $\psi \oplus \lambda$.
This is an equivalence relation on the subset of semi-invertible
extensions in $\Hom_G(A,Q(B))$ and the corresponding equivalence
classes form an abelian group which we denote by $\Ext^{ -1/2}(A,B)$.
For any locally compact space $X$ we consider $C_0(X) \otimes A$ as a
$G$-algebra with the trivial $G$-action on the tensor factor $C_0(X)$.
When $X = \ (0,1]$ we denote $C_0(0,1] \otimes A$ by $\cone(A)$.
Similarly, we set $SA = C_0(0,1)\otimes A$.

\begin{lemma}\label{crux} Let $\lambda : \cone(A) \to Q(B)$ be a
$G$-extension. It follows that there is an asymptotic homomorphism
$\pi = (\pi_t)_{t \in [1,\infty)} : \cone(A) \to M_2(M(B))$ such that
$$
q_{M_2(B)} \circ \pi_t = \left ( \begin{smallmatrix} \lambda &  {}  \\
{}  & 0 \end{smallmatrix} \right )
$$
for all $t \in [1,\infty)$.
\begin{proof} The proof is based on an idea of Voiculescu, cf. \cite{V}.
Let $\mu : \cone(A) \to M(B)$ be a continuous, self-adjoint and
homogeneous lift of $\lambda$ such that $\|\mu(x)\| \leq 2 \|x\|$ for
all $x \in \cone(A)$. Such $\mu$ exists by the Bartle-Graves selection
theorem, cf. \cite{L}. Define $\varphi_{s} : \cone(A) \to \cone(A)$ such
that $\varphi_{s}(f)(t) = f((1-s)t),  s \in [0,1]$. Choose continuous
functions $f_i : [1, \infty) \to [0,1], i = 0,1,2,\cdots$, such that
\begin{enumerate}
\item[1)] $f_0(t) =0$ for all $t \in [1,\infty)$,
\item[2)] $f_n \leq f_{n+1}$ for all $n$,
\item[3)] for each $ n \in \Bbb N$, there is an $m_n \in \Bbb N$ such
that $f_i(t) = 1$ for all $i \geq m_n$, and all $t \in [1,n+1]$,
\item[4)] $\lim_{t \to \infty} \max_i |f_i(t) - f_{i+1}(t)| = 0$.
\end{enumerate}
Let $F_1 \subseteq F_2 \subseteq F_3 \subseteq \cdots $ be an
increasing sequence of finite subsets with dense union in $\cone(A)$.
Write $G = \bigcup_n K_n$ where $K_1 \subseteq K_2 \subseteq K_3
\subseteq \cdots $ are compact subsets of $G$. For each $n$, choose
$m_n \in \Bbb N$ as in 3). We may assume that $m_{n+1} > m_n$.
By Lemma 1.4 of \cite{K} we can choose elements
$$
X^n_{0}   \geq X^n_{1} \geq X^n_{2} \geq  \cdots
$$
in $B$ such that $0 \leq X^n_i \leq 1$ for all $i$ and $X^n_i = 0$
for $i \geq  m_n$, and
\begin{enumerate}
\item[1')] $X^n_{i}X^n_{i+1} = X^n_{i+1}$ for all $i$,
\item[2')] $\|X^n_i b - b\| \leq \frac{1}{n}$ for all $i = 0,1,2,
\cdots , m_n - 1$, and all $b \in S_n$,
\item[3')] $\|X^n_i y - y X^n_i\| \leq \frac{1}{n}$ for all $ i$ and
all $y \in L_n$,
\item[4')] $\|g \cdot X_i^n - X_i^n\| \leq \frac{1}{n}, \ g \in K_n$,
for all $i$,
\item[5')] $\|X_i^n ( g \cdot \mu(a) - \mu(g \cdot a)) -
( g \cdot \mu(a) - \mu(g \cdot a))\| \leq \frac{1}{n}, \ g \in K_n,
a \in F_n$, for all $i = 0,1,2, \cdots , m_n -1$,
\end{enumerate}
where $L_n$ and $S_n$ are the compact sets $L_n = \{ \mu(\varphi_{s}(a)) :
\ s \in [0,1], \ a \in F_n \}$ and
\eQ{\spl{
&S_n = \{\mu(\varphi_{s}(a)) + \mu(\varphi_{s}(b)) -
\mu(\varphi_{s}(a + b)) : \ a,b \in F_n, \ s \in [0,1] \} \\
& \cup \{\mu(\varphi_{s}(ab)) - \mu(\varphi_{s}(a))\mu(\varphi_{s}(b)) :
\ a,b \in F_n, \ s \in [0,1] \} \ . }}
Since we choose the $X$'s recursively we can arrange that
$X^{n+1}_i X^n_k = X^n_k$ for all $k$ and all $i \leq m_{n+1}$.
By connecting first $X^n_0$ to $X^{n+1}_0$ via the straight line
between them, then $X^n_{1}$ to $X^{n+1}_{1}$ via a straight line,
then $X^n_{2}$ to $X^{n+1}_{2}$ etc., we obtain norm-continuous
pathes, $X(t,i), t \in [n,n+1], i = 0,1,2,3, \cdots$, in $B$ such
$X(n,i) = X^n_i, \ X(n+1,i) = X^{n+1}_i$ for all $i$ and
\begin{enumerate}
\item[a)] $X(t,{i})X(t,i+1) = X(t,i+1), \  t \in [n,n+1]$, for all $i$,
\item[b)] $\|X(t,i) b - b\| \leq \frac{1}{n}$ for all $i = 0,1,2, \cdots ,
m_n - 1, \ t \in [n,n+1]$ and all $b \in S_n$,
\item[c)] $\|X(t,i) y - y X(t,i)\| \leq \frac{1}{n}$ for all $ i $,
all $t \in [n,n+1]$ and all $y \in L_n$,
\item[d)] $\| g \cdot X(t,i) - X(t,i)\| \leq \frac{1}{n},
g \in K_n, t \in [n,n+1]$, for all $i$,
\item[e)] $\|X(t,i)( g \cdot \mu(a) - \mu(g \cdot a)) -
( g \cdot \mu(a) - \mu(g \cdot a))\| \leq \frac{1}{n},
\ g \in K_n, a \in F_n, t \in [n,n+1]$, for all $i = 0,1,\cdots , m_n-1$.
\end{enumerate}
In addition, $X(t,i) = 0, i \geq m_{n+1}, t \in [n,n+1]$. Let $l_2(B)$
denote the Hilbert $B$-module of sequences $(b_1,b_2,b_3, \cdots )$ in
$B$ such that $\sum_{i=1}^{\infty} b_i^*b_i$ converges in norm. Writing
an element $(b_1,b_2,b_3, \cdots ) \in l_2(B)$ as the sum
$\sum_{i=0}^{\infty} b_ie_i$ we define a representation $V$
of $G$ on $l_2(B)$ such that $V_g(\sum_{i=0}^{\infty} b_ie_i) =
\sum_{i=0}^{\infty} (g \cdot b_i)e_i$. Then $G$ acts by automorphisms
on $\Bbb L(l_2(B))$ ( = the adjoinable operators on $l_2(B)$) such that
$g \cdot m = V_gmV_{g^{-1}}$.  Set
$$
T_t = \left ( \begin{matrix} \sqrt{1 - X(t,0)} &  \sqrt{X(t,0) -
X(t,1)} &  \sqrt{X(t,1) - X(t,2)} & \hdots \\ 0 & 0 & 0 & \hdots \\
0 & 0 & 0 & \hdots \\ \vdots & \vdots & \vdots & \ddots \end{matrix}
\right ) \in \Bbb L(l_2(B)) .
$$
Then $P_t = T_t^*T_t$ is a projection in $\Bbb L(l_2(B))$ since
$T_tT_t^*$ clearly is. Note that $P_t$ is tri-diagonal because of
condition a) above, and that the entries of $P_t$ are all in $B$,
with the notable exception of the $1 \times 1$-entry which is equal
to $1$ modulo $B$. We define $\delta_t : \cone(A) \to \Bbb L(l_2(B))$ by
$$
\delta_t(a)(\sum_{i=0}^{\infty} b_ie_i) = \sum_{i=0}^{\infty}
\mu(\varphi_{f_i(t)}(a))b_i e_i  .
$$
Set $\pi_t(a) = P_t \delta_t(a)P_t$ for $a \in \cone(A)$ and
$t \in [1,\infty)$. We assert that $\pi = (\pi_t)_{t \in [1,\infty)}$
is an asymptotic homomorphism. By using the continuity of $\mu$ and
that $\{\varphi_s(a) : s \in [0,1]\}$ is a compact set for fixed $a$,
it follows readily that the family of maps $a \mapsto \pi_t(a),
t \in [1,\infty)$, is an equicontinuous family. Since each $\pi_t$
is self-adjoint and homogeneous, it suffices therefore to take an $n$
and elements $a,b \in F_n$, $g \in K_n$, and check that
$$
\lim_{t \to \infty} P_t \delta_t(a)P_t\delta_t(b)P_t -
P_t \delta_t(ab)P_t = 0,
$$
$$
\lim_{t \to \infty} P_t\delta_t(a + b)P_t - P_t\delta_t(a)P_t -
P_t\delta_t(b)P_t = 0 ,
$$
and
$$
\lim_{t \to \infty} P_t \delta_t ( g \cdot a) P_t - g \cdot
(P_t \delta_t(a)P_t) = 0.
$$
The first two limits are zero by 4), b) and c), the third by d) and e).
For each $a,t$, $P_t\delta_t(a)P_t = \diag ( \mu(a), 0,0, \cdots)$
modulo $\Bbb K(l_2(B))$ ( = the ideal of 'compact' operators on $l_2(B)$).
Since $B$ is weakly stable there is an equivariant isomorphism
$l_2(B) \simeq B \oplus B$ of Hilbert $B$-modules which leaves
the first coordinate invariant. We can therefore transfer $\pi$
to an asymptotic homomorphism $\pi = (\pi_t)_{t \in [1,\infty)} :
\cone(A) \to \Bbb L(B \oplus B) = M_2(M(B))$ with the stated property.

\end{proof}
\end{lemma}

Two $G$-extensions $\varphi, \psi \in \Hom_G(A,Q(B))$ are \it
strongly homotopic \rm when there is a path $\bold\Phi_t \in
\Hom_G(A,Q(B)), t \in [0,1]$, such that $\bold\Phi_0 = \varphi,
\bold\Phi_1 = \psi$ and $t \mapsto \bold\Phi_t(a)$ is continuous for
all $a \in A$.

\begin{thm}\label{cor} Let $\varphi : A \to Q(B)$ be a $G$-extension
which is strongly homotopic to $0$ in $\Hom_G (A, Q(B))$. It follows
that there is an asymptotic homomorphism
$\pi = (\pi_t)_{t \in [1,\infty)} : A \to M_2(M(B))$ such that
$$
q_{M_2(B)} \circ \pi_t = \left ( \begin{smallmatrix} \varphi &  {}  \\
{}  & 0 \end{smallmatrix} \right )
$$
for all $t \in [1,\infty)$.
\begin{proof} Since $\varphi$ is strongly homotopic to $0$ there is an
equivariant $*$-homomorphism $\mu : A \to \cone(D)$, where
$D \subseteq Q(B)$ is a separable $G$-algebra containing $\varphi(A)$,
and an equivariant $*$-homomorphism $\lambda : \cone(D) \to Q(B)$ such
that $\varphi = \lambda \circ \mu$. Apply Lemma \ref{crux} to $\lambda$.
\end{proof}
\end{thm}

\begin{cor}\label{A!} Every $G$-extension $\varphi : SA \to Q(B)$ is
semi-invertible.
\begin{proof} Let $\alpha \in \Aut SA$ be the automorphism of $SA$ given
by $\alpha(f)(t) = f(1-t)$. It is wellknown that
$\varphi \oplus (\varphi \circ \alpha)$ is strongly homotopic to $0$.
Hence $\varphi \oplus (\varphi \circ \alpha) \oplus 0$ is asymptotically
split by Theorem \ref{cor}.
\end{proof}
\end{cor}

Because of Corollary \ref{A!} we drop the superscript $-1/2$ and write
$\Ext(SA,B)$ instead of $\Ext^{-1/2}(SA,B)$.

\begin{lemma}\label{LL1} Let $\varphi, \psi : SA \to Q(B)$ be two
$G$-extensions which are strongly homotopic. It follows that $\varphi$
and $\psi$ are stably unitarily equivalent.
\begin{proof} It follows from Theorem \ref{cor} that
$\lambda_1 = (\varphi \circ \alpha) \oplus  \varphi \oplus 0$ and
$\lambda_2 = (\varphi \circ \alpha) \oplus \psi \oplus 0$ are both
asymptotically split. Since $\psi \oplus \lambda_1$ and
$\varphi \oplus \lambda_2$ are unitarily equivalent, the conclusion
follows because infinite direct sums are well-defined for asymptotically
split extensions.
\end{proof}
\end{lemma}

Set $IB = C[0,1] \otimes B$ and let $e_t : IB \to B$ denote evaluation
at $t \in [0,1]$ and note that $e_t$ defines a equivariant
$*$-homomorphisms $M(IB) \to M(B)$ and $Q(IB) \to Q(B)$ which we again
denote by $e_t$. Two $G$-extensions $\varphi, \psi \in \Hom_G(A,Q(B))$
are \it homotopic \rm when there is a $G$-extension
$\bold\Phi \in \Hom_G(A,Q(IB))$ such that $e_0 \circ \bold\Phi =
\varphi$ and $e_1 \circ \bold\Phi = \psi$. As in \cite{MT2} we denote
the set of homotopy classes of $G$-extensions by $\Ext(A,B)_h$. In general
this is merely an abelian semigroup, but $\Ext(SA,B)_h$ is a group.

The Connes-Higson construction associates to any $G$-extension
$\varphi \in \Hom_G(A,Q(B))$ an asymptotic homomorphism $CH(\varphi) :
SA \to B$ in the following way, cf. \cite{CH}, \cite{GHT}: By Lemma 1.4
of \cite{K} or Lemma 5.3 of \cite{GHT} there is a norm-continuous path
$\{u_t\}_{t \in [1,\infty)}$ of elements in $B$ such that
$0 \leq u_t \leq 1$ for all $t$, $\lim_{t \to \infty} \|u_t b - b\| = 0$
for all $b \in B$, $\lim_{t \to \infty} \|u_tm - m u_t \| = 0$ for all
$m \in q_B^{-1}(\varphi(A))$ and $\lim_{t \to \infty} \|g \cdot u_t - u_t\|
= 0$ for all $g \in G$. From these data $CH(\varphi)$ is determined up
to asymptotic equality as the equicontinuous\footnote{Equicontinuity of
an asymptotic homomorphism $\pi = (\pi_t)_{t \in [1,\infty)} : A \to B$
means that $A \times G \ni (a,g) \mapsto g \cdot \pi_t(a), t \in
[1,\infty)$, is an equicontinuous family of maps.} asymptotic homomorphism
$CH(\varphi) : SA \to B$ which satisfies that
$$
\lim_{t \to \infty} CH(\varphi)_t ( f \otimes a) - f(u_t)x  = 0, \ \ \ \
x \in q_B^{-1}(\varphi(a)),
$$
for all $f \in C_0(0,1)$ and all $a \in A$. Let $[[SA,B]]$ denote the
abelian group of homotopy classes of asymptotic homomorphisms, $SA \to B$,
cf. \cite{CH}, \cite{GHT}. The Connes-Higson construction defines in the
obvious way a semi-group homomorphism $CH : \Ext(A,B)_h \to [[SA,B]]$.
Since there is a canonical (semi-group) homomorphism $\Ext^{-1/2}(A,B)
\to \Ext(A,B)_h$ we may also consider the Connes-Higson construction
as a homomorphism $CH :\Ext^{-1/2}(A,B) \to [[SA,B]]$. Notice that
$\Ext(SA,B)$ and $\Ext(SA,B)_h$ are both abelian groups and the canonical
map $\Ext(SA,B) \to \Ext(SA,B)_h$ is a surjective group homomorphism by
Corollary \ref{A!}. In Corollary \ref{coinc} below we show that it is an
isomorphism.

\section{On equivalence of asymptotic homomorphisms}

\begin{lemma}\label{!} Let $A$ and $B$ be separable $G$-algebras, $B$
weakly stable. Let $\varphi = (\varphi_t)_{t \in [1,\infty )} : A \to B$
be an asymptotic homomorphism which is homotopic to $0$. It follows that
there is an asymptotic homomorphism $\psi = (\psi_t)_{t \in [1,\infty)}
: A \to B$ and a norm-continuous path $\{W_t\}_{t \in [1,\infty)}$ of
$G$-invariant unitaries in $M(M_2(B))$ such that
$$
\lim_{t \to \infty} \left ( \begin{smallmatrix} \varphi_t(a) & {}  \\
{}  & \psi_t(a) \end{smallmatrix} \right ) -
W_t\left ( \begin{smallmatrix} 0 & {}  \\ {}  & \psi_t(a)
\end{smallmatrix} \right )W_t^* = 0
$$
for all $a \in A$.
\begin{proof} Let $\bold\Phi = (\bold\Phi_t)_{t \in [1,\infty)} :
A \to IB$ be an asymptotic homomorphism such that $e_0 \circ
\bold\Phi_t(a) = 0,  e_1 \circ \bold\Phi_t(a) = \varphi_t(a)$ for all
$t \in [1,\infty),  a \in A$. We may assume that both $\varphi$ and
$\bold\Phi$ are equicontinuous, cf. Proposition 2.4 of \cite{Th2}.
Let $F_1 \subseteq F_2 \subseteq F_3 \subseteq \cdots $ be a sequence
of finite subsets with dense union in $A$. For each $n$ there is
$\delta_n > 0$ with the property that
$$
\|e_x \circ \bold\Phi_t(a) - e_y \circ \bold\Phi_t(a)\| < \frac{1}{n}
$$
when $|x -y | < \delta_n,  t \in [1,n],  a \in F_n$. Choose then a
sequence of functions $f_k : [1,\infty) \to [0,1]$ such that $f_1(t) =1,
\ f_k  \geq f_{k+1}, \ |f_k(t) - f_{k+1}(t)| < \delta_n, \ t \in [1,n]$
for all $k,n$ and such that $f_k|_{[1,n]} = 0$ for all but finitely many
$k$'s for all $n$. Set $\lambda^n_t(a) = e_{f_n(t)} \circ \bold\Phi_t(a)$
for all $a \in A,  n \in \Bbb N, t \in [1,\infty)$. Note
that $\|\lambda^i_t(a) - \lambda^{i+1}_t(a)\| < \frac{1}{n}, a \in F_n,
\ t \in [1,n]$, for all $i$ and $n$.
Then
$$
\mu_t(a) = \diag (\varphi_t(a), \lambda^1_t(a),\lambda^2_t(a),
\lambda^3_t(a), \cdots ) \ \in \ \Bbb K(l_2(B))
$$
and
$$
\delta_t(a) = \diag (0, \lambda^1_t(a),\lambda^2_t(a), \lambda^3_t(a),
\cdots ) \ \in \ \Bbb K(l_2(B))
$$
define asymptotic homomorphisms $\mu , \delta : A \to  \Bbb K(l_2(B))$.
By connecting appropriate permutation unitaries, acting on $l_2(B)$ by
permutations of $B$-coordinates, we get a norm-continuous path of
$G$-invariant unitaries $\{S_t\}_{t \in [1,\infty)} \subseteq
\Bbb L (l_2(B))$ such that
$$
S_t \delta_t(a) S_t^* =  \diag (\lambda^1_t(a),\lambda^2_t(a),
\lambda^3_t(a), \cdots )
$$
for all $a,t$. Then $\lim_{t  \to \infty} \mu_t(a) -
S_t\delta_t(a)S_t^* = 0$ for all $a \in A$. Since $B$ is weakly
stable there is an isomorphism $l_2(B) \to B \oplus B$ of Hilbert
$B,G$-algebras which fixes the first coordinate. Applying this
isomorphism in the obvious way and remembering the identifications
$\Bbb K(B \oplus B) = M_2(B)$ and $\Bbb L (B \oplus B) = M(M_2(B))$
gives the result.
\end{proof}
\end{lemma}

\begin{thm}\label{!!} Let $A$ and $B$ be separable $G$-algebras,
$B$ weakly stable. Assume that $[[A,B]]$ is a group. Two asymptotic
homomorphisms, $\varphi = (\varphi_t)_{t \in [1,\infty)}, \ \psi =
(\psi_t)_{t \in [1,\infty)} : A \to B$, are homotopic if and only if
there is an asymptotic homomorphism $\lambda = (\lambda_t)_{t \in
[1,\infty)} : A \to B$ and a norm-continuous path $\{W_t\}_{t \in
[1,\infty)}$ of $G$-invariant unitaries in $M(M_2(B))$ such that
$$
\lim_{t \to \infty} \left ( \begin{smallmatrix} \varphi_t(a) & {}  \\
{}  & \lambda_t(a) \end{smallmatrix} \right ) -
W_t\left ( \begin{smallmatrix} \psi_t(a) & {}  \\ {}  &
\lambda_t(a) \end{smallmatrix} \right )W_t^* = 0
$$
for all $a \in A$.
\begin{proof} The 'if' part is easy and the 'only if' part follows
from Lemma \ref{!} in the same way as Lemma \ref{LL1} follows from
Theorem \ref{cor}.
\end{proof}
\end{thm}

\begin{lemma}\label{L2!} Let $B$ be a weakly stable $G$-algebra and
$D_0$ a separable $G$-subalgebra of $C_b([1,\infty),B)$. Let $V_1,V_2,
\cdots, V_N \in M(B)$ be $G$-invariant isometries. There is then a weakly
stable separable $G$-subalgebra $D$ of $C_b([1,\infty),B)$ such that
$V_iD \cup V_i^*D \cup D_0 \subseteq D$ for all $i = 1,2, \cdots , N$.
\begin{proof} Since $B$ is weakly stable we can write
$B = B \otimes \mathcal K$ with $G$ acting trivially on the
tensor-factor $\mathcal K$. We embed $\mathcal K$ into
$M(B \otimes \mathcal K)$ via $x \mapsto 1_B \otimes x$. Let
$\{f_n\} \subseteq C_b([1,\infty),B \otimes \mathcal K)$ be a dense
sequence in $D_0$. For each $n \in \Bbb N$ there is a function
$g_n \in C_b([1,\infty),\mathcal K)$ such that $\|g_nf_n -f_n\| <
\frac{1}{n}$. Let $E_{00}$ be the $C^*$-algebra generated by
$\{g_n\}_{n=1}^{\infty}$. Then $E_{00} \subseteq C_b([1,\infty),\mathcal K)
\subseteq C_b([1,\infty),B^+ \otimes \mathcal K)$. Consider a positive
element $f \in E_{00}$ and an $\epsilon > 0$. Set $U_j = ]j,j+2[ \cap
[1,\infty[, j = 0,1,2, \cdots$. We can then find a sequence
$p_0 \leq p_1 \leq p_2 \leq \cdots$ of projections in $\mathcal K$
such that
$$
\sup_{x \in \overline{U_j}} \|p_jf(x)p_j - f(x)\| < \epsilon  .
$$
Let $\{h_j\}$ be a partition of unity in $C_b[1,\infty)$ subordinate
to the cover $\{U_j\}$ and set $g(t) = \sum_{j=0}^{\infty} h_j(t)
p_jf(t)p_j$. Then $g \in C_b([1,\infty),\mathcal K)$, $g \geq 0,
\|g - f\| < \epsilon$. For each $j$ we choose a partial isometry
$v_j \in \mathcal K$ such that $v_jv_j^* = p_{j+2}$, $v_j^*v_jp_{j+2} = 0$
and $v_j^*v_jv_k^*v_k = 0,  k < j$. Set $h(t) = \sum_{j=0}^{\infty}
\sqrt{h_j(t)}v_j$. Then $hh^*g = g$ and $h^*hg = 0$. It follows that we
can find a sequence $E_{00} = X_1 \subseteq X_2 \subseteq X_3 \subseteq
\cdots $ of separable $C^*$-subalgebras of $C_b([1,\infty),\mathcal K)$
and for each $n$ have a dense sequence $\{f_1,f_2, \cdots \}$ in the
positive part of $X_n$ and elements $\{v_1,v_2, \cdots \}$ in $X_{n+1}$
such that $\|f_k - v_k^*v_k\| < \frac{1}{k}$ and $v_k^*v_kv_kv_k^* = 0$
for all $k$. It follows then from Proposition 2.2 and Theorem 2.1 of
\cite{HR} that $E_0 = \overline{\bigcup_n X_n}$ is a separable stable
$C^*$-subalgebra of $C_b([1,\infty),\mathcal K)$ such that $E_{00}
\subseteq E_0$. Note that $E_0$ contains a sequence $\{r_n\}$ with
the property that $\lim_{n \to \infty} r_n x = x$ for all $x \in D_0$
since $E_{00}$ does. Set $W = \{V_1,V_2, \cdots , V_N\}
\cup \{V_1^*,V_2^*, \cdots , V_N^*\}$. By repeating the above argument
with $D_0$ substituted by the $G$-algebra $D_1$ generated by $D_0 \cup
WD_0 \cup E_0D_0$, we get a stable $C^*$-subalgebra $E_1 \subseteq
C_b([1,\infty),\mathcal K)$ which contains a sequence $\{r_n\}$ such
that $\lim_{n \to \infty}r_ny = y$ for all $y \in D_1$. It is clear
from the construction that we can arrange that $E_0 \subseteq E_1$.
We can therefore continue this procedure to obtain sequences of separable
$G$-algebras, $D_0 \subseteq D_1 \subseteq D_2 \subseteq D_3 \subseteq
\cdots$ in $C_b([1,\infty), B \otimes \mathcal K)$, and $E_0 \subseteq
E_1 \subseteq E_2 \subseteq E_3 \subseteq \cdots $ in $C_b([1,\infty),
\mathcal K) \subseteq C_b([1,\infty),B^+\otimes \mathcal K)$ such that
each $E_n$ is stable and contains a sequence $\{r_k\}$ such that
$\lim_{k \to \infty} r_k x = x, \ x \in D_n$, and $D_n \cup WD_n \cup
E_nD_n \subseteq D_{n+1}$ for all $n$. Set
$E_{\infty} = \overline{ \bigcup_n E_n}$ and $D = \overline{ \bigcup_n D_n}$.
It follows from Corollary 4.1 of \cite{HR} that $E_{\infty}$ is stable. By
construction $V_iD \cup V_i^*D \subseteq D$ for all $i$ and $E_{\infty}D
\subseteq D$. The last property ensures that $D$ is an ideal in the
$G$-algebra $E$ generated by $E_{\infty}$ and $D$. There is therefore a
$*$-homomorphism $ \lambda : E_{\infty} \to M(D)$. By construction an
approximate unit for $E_{\infty}$ is also an approximate unit for $D$ so
$\lambda$ extends to a $*$-homomorphism $\lambda : M(E_{\infty}) \to M(D)$
which is strictly continuous on the unit ball of $M(E_{\infty})$. Since
$E_{\infty}$ is stable there is a sequence $P_i, i = 1,2, \cdots $, of
orthogonal and Murray-von Neumann equivalent projections in $M(E_{\infty})$
which sum to $1$ in the strict topology. Then $Q_i = \lambda(P_i),
i = 1,2,\cdots $, is a sequence of orthogonal and Murray-von Neumann
equivalent projections in $M(D)$ which sum to $1$ in the strict topology.
Since $E_{\infty}$ consists entirely of $G$-invariant elements it follows
that all the $Q_i$'s are $G$-invariant. Consequently $D \simeq Q_1DQ_1
\otimes \mathcal K$ as $G$-algebras, proving that $D$ is weakly stable.
\end{proof}
\end{lemma}

Two asymptotic homomorphisms $\varphi = (\varphi_t)_{t \in [1,\infty)}, \
\psi = (\psi_t)_{t \in [1,\infty)} : A \to B$ will be called \it
equi-homotopic \rm when there is a family $\bold\Phi^{\lambda} =
(\bold\Phi^{\lambda}_t)_{t \in [1,\infty)} : A \to B, \ \lambda \in
[0,1]$, of asymptotic homomorphisms such that the family of maps,
$[0,1] \ni \lambda \mapsto \bold\Phi_t^{\lambda}(a), \ t \in [1,\infty)$,
is equicontinuous for each $a \in A$.

\begin{thm}\label{TT1} Let $A$ and $B$ be separable $G$-algebras,
$B$ weakly stable. Let $\varphi = (\varphi_t)_{t \in [1,\infty)}, \
\psi = (\psi_t)_{t \in [1,\infty)} : SA \to B$ be asymptotic homomorphisms.
Then the following are equivalent:
\begin{enumerate}
\item[1)] $\varphi$ and $\psi$ are homotopic $($i.e. $[\varphi] = [\psi]$
in $[[SA,B]]$$)$.
\item[2)] $\varphi$ and $\psi$ are equi-homotopic.
\item[3)] There is an asymptotic homomorphism $\lambda =
(\lambda_t)_{t \in [1,\infty)} : SA \to B$ and a norm-continuous
path $\{W_t\}_{t \in [1,\infty)}$ of $G$-invariant unitaries in
$M(M_2(B))$ such that
$$
\lim_{t \to \infty} \left ( \begin{smallmatrix} \varphi_t(a) & {}  \\
{}  & \lambda_t(a) \end{smallmatrix} \right ) -
W_t\left ( \begin{smallmatrix} \psi_t(a) & {}  \\ {}  & \lambda_t(a)
\end{smallmatrix} \right )W_t^* = 0
$$
for all $a \in A$.
\end{enumerate}
\begin{proof} The equivalence 1) $\Leftrightarrow$ 3) follows from
Theorem \ref{!!} and the implication 2) $\Rightarrow$ 1) is trivial,
so we need only prove that 1) $\Rightarrow$ 2). To this end, let
$[[SA,B]]^e$ denote the set of equi-homotopy classes of asymptotic
homomorphisms $SA \to B$. Choose $G$-invariant isometries $V_1,V_2
\in M(B)$ such that $V_1V_1^* + V_2V_2^* = 1$ and define a composition
in $[[SA,B]]^e$ by
$$
[\varphi] + [\psi] = [(V_1\varphi_tV_1^* + V_2\psi_tV_2^*)_{t \in
[1,\infty)}] \ .
$$
It follows from Lemma \ref{L2!} that $[[SA,B]]^e$ is a group. It
suffices therefore to show that the natural map $[[SA,B]]^e \to
[[SA,B]]$ has trivial kernel. If $\varphi$ is an asymptotic homomorphism
representing an element in the kernel we conclude from Lemma \ref{!}
that there is a norm-continuous path $W_t, t \in [1,\infty)$, of
$G$-invariant unitaries in $M_2(M(B)))$ and an asymptotic homomorphism
$\psi$ such that
$$
\lim_{t \to \infty} \left ( \begin{smallmatrix} \varphi_t(a) & {} & {} &
{} \\ {}  & \psi_t(a) & {} & {} \\ {} & {} & 0 & {} \\ {} & {} & {} & 0
\end{smallmatrix} \right ) - \left ( \begin{smallmatrix} W_t & {} \\
{} & W_t^* \end{smallmatrix} \right ) \left ( \begin{smallmatrix} 0 &
{} & {} & {} \\ {}  & \psi_t(a) & {} & {} \\ {} & {} & 0 & {} \\ {} &
{} & {} & 0 \end{smallmatrix} \right ) \left ( \begin{smallmatrix}
W_t^* & {} \\ {} & W_t \end{smallmatrix} \right )    = 0
$$
for all $a \in SA$. By a standard rotation argument we can remove the
unitaries $\left ( \begin{smallmatrix} W_t & {} \\ {} & W_t^*
\end{smallmatrix} \right )$ via an equi-homotopy and we see in this way
that $[\varphi] + [\psi] = [\psi]$ in $[[SA,B]]^e$. Hence $[\varphi] = 0$
in $[[SA,B]]^e$.

\end{proof}
\end{thm}

Simple examples show that the implications 1) $\Rightarrow$ 2) and
1) $\Rightarrow$ 3) of Theorem \ref{TT1} generally fail in $[[A,B]]$.

\section{Making genuine homomorphisms out of asymptotic ones}

Let $A$ and $B$ be separable $C^*$-algebras. Set
$$
M(B)_G = \{ x \in M(B) : \ G \ni g \mapsto g \cdot x \
\text{is norm-continuous} \}
$$
and
$$
Q(B)_G = \{ x \in Q(B) : \ G \ni g \mapsto g \cdot x \
\text{is norm-continuous} \} .
$$
Then
\eq{\label{ex1}\begin{xymatrix}{
0  \ar[r]  &  B  \ar[r]  &  M(B)_G \ar[r]  &  Q(B)_G \ar[r]  & 0 \  }
\end{xymatrix} }
is a short exact sequence of $G$-algebras. (This is not trivial -
the surjectivity of the quotient map follows from Theorem 2.1 of
\cite{Th1}.) We are going to construct a map $\alpha :
[[SA,Q(B)_G \otimes \mathcal K]] \to \Ext(SA,B \otimes \mathcal K)_h$.
The key to this is another variant of the Voiculescu's tri-diagonal
projection trick from \cite{V}. Let $b$ be a strictly positive element
of $B \otimes \mathcal K$, $0 \leq b \leq 1$. A \it unit sequence \rm
in $B \otimes \mathcal K$ is a sequence $\{u_n\}_{n=0}^{\infty} \subseteq
B \otimes \mathcal K$ such that
\begin{enumerate}
\item[0)] there is a continuous function $f_n : [0,1] \to [0,1]$ which
is zero in a neighbourhood of $0$ and $u_n = f_n(b)$,
\item[1)] $0 \leq u_n \leq 1$ for all $n = 0,1,2,3, \cdots $,
\item[2)] $u_{n+1}u_n = u_n $ for all $n$,
\item[3)] $\lim_{n \to \infty} u_n x = x, \ x \in B \otimes \mathcal K$,
\item[4)] $\lim_{n \to \infty} \|g \cdot u_n - u_n\| = 0, \ g \in G$.
\end{enumerate}

Let $\{e_{ij}\}_{i,j = 0}^{\infty}$ be the matrix units acting on
$l_2(B \otimes \mathcal K)$ in the standard way.

\begin{lemma}\label{X1} Let $\mathcal U = \{u_n\}$ be a unit sequence
in $B \otimes \mathcal K$. Then
$$
\sqrt{u_0} e_{00} + \sum_{j=1}^{\infty} \sqrt{u_j - u_{j-1}}e_{0 j}
$$
converges in the strict topology to a partial isometry $V$ in
$\Bbb L(l_2(B \otimes \mathcal K))$ such that $VV^* = e_{00}$.
\begin{proof} Let $b = (b_0,b_1,b_2, \cdots ) = \sum_{i = 0}^{\infty}
b_ie_i \in l_2(B \otimes \mathcal K)$. Then
\eQ{\spl{
&\|\sum_{j =n}^m \sqrt{u_j - u_{j-1}} e_{0j}(b)\|^2 =
\| \sum_{k,j = n}^m b_k^* \sqrt{u_k - u_{k-1}} \sqrt{u_j - u_{j-1}} b_j  \|
\\
& \\
& = \|\sum_{k =n}^m b_k^* (u_k - u_{k-1}) b_k +  \sum_{ k =n}^{m-1}
b_k^*\sqrt{u_k - u_{k-1}} \sqrt{u_{k+1} - u_{k}} b_{k+1}  + \\
& \   \ \ \ \ \ \ \ \ \ \ \ \ \ \ \ \sum_{ k =n}^{m-1} b_{k+1}^*
\sqrt{u_{k+1} - u_{k}} \sqrt{u_{k} - u_{k-1}} b_{k}\| \\
& \\
& \leq \|\sum_{k=n}^m b_k^*b_k \| +  2 \sqrt{ \|\sum_{k =n}^{m-1}
b_k^*b_k\|}  \sqrt{ \|\sum_{k =n +1 }^{m} b_k^*b_k\|} , }}
proving that $\sum_{j=1}^{\infty}  \sqrt{u_j - u_{j-1}}e_{0 j}(b)$
converges in $l_2(B \otimes \mathcal K)$. And
\eQ{\|(\sum_{j =n}^m \sqrt{u_j - u_{j-1}} e_{0j})^*(b)\|^2 =
\| \sum_{j=n}^m b_0^*(u_j - u_{j-1})b_0\| , }
proving that also $(\sum_{j=1}^{\infty}  \sqrt{u_j -
u_{j-1}}e_{0 j})^*(b)$ converges in $l_2(B \otimes \mathcal K)$.
It follows that
$$
V =
\sqrt{u_0} e_{00} + \sum_{j=1}^{\infty} \sqrt{u_j - u_{j-1}}e_{0 j}
$$
exists as a strict limit in $\Bbb L(l_2(B \otimes \mathcal K))$.
It it then straightforward to check that $VV^* = e_{00}$.
\end{proof}
\end{lemma}

Let $P_{\mathcal U} = V^*V$ and note that $P_{\mathcal U}$ is
tri-diagonal with respect to the matrix units $\{e_{ij}\}$. Fix now a
continuous and homogeneous section $\chi$ for the map $q_B \otimes
\id_{\mathcal K} : M(B)_G \otimes \mathcal K \to Q(B)_G \otimes
\mathcal K$. Consider an equicontinuous asymptotic homomorphism
$\varphi = (\varphi_t)_{t \in [1,\infty)} : A \to Q(B)_G \otimes
\mathcal K$. Let $F_1 \subseteq F_2 \subseteq F_3 \subseteq \cdots $
be a sequence of finite sets with dense union in $A$ and $K_1 \subseteq
K_2 \subseteq K_3 \subseteq \cdots$ a sequence of compact subsets in $G$
such that $\bigcup_n K_n = G$. It is easy to see that there is a unit
sequence $\{u_n\}$ in $B \otimes \mathcal K$ with the following properties :
\begin{enumerate}
\item[5)] $ \|u_n   \chi(\varphi_t(a)) -  \chi(\varphi_t(a))u_n\| \leq
\frac{1}{n}, \ a \in F_n, \ t \in [1,n +1]$,
\item[6)] $\|(1 - u_{n})( \chi(\varphi_t(ab)) - \chi(\varphi_t(a))
\chi(\varphi_t(b)))\| \leq \|\varphi_t(ab) - \varphi_t(a)\varphi_t(b)\|
+ \frac{1}{n}, \ t \in [1,n+1], \ a,b \in F_n$,
\item[7)] $\| (1 - u_{n})(\chi(\varphi_t(a + b)) - \chi(\varphi_t(a)) -
\chi(\varphi_t(b)))\| \leq \|\varphi_t(a + b) - \varphi_t(a) -
\varphi_t(b)\| + \frac{1}{n},  \ t \in [1,n+1], \ a,b \in F_n $,
\item[8)] $\|( 1- u_n)(g \cdot \chi(\varphi_t(a)) -
\chi(\varphi_t(g \cdot a)))\| \leq \| g \cdot \varphi_t(a) -
\varphi_t(g \cdot a)\| + \frac{1}{n}, \ t \in [1,n], \ a \in F_n,
\ g \in K_n$.

\end{enumerate}

Let $\{\varphi_{t_n}\}_{n \in \Bbb N}$ be a discretization of
$\varphi$, cf. Lemma 5.1 of \cite{MT1}, such that
\begin{enumerate}
\item[9)] $t_n \leq n$ for all $n \in \Bbb N$.
\end{enumerate}
Set
$$
\widetilde{\varphi}(a) = P_{\mathcal U} ( \sum_{j = 0}^{\infty}
\chi(\varphi_{t_{j+1}}(a)) e_{jj})P_{\mathcal U} \ .
$$
Then $\widetilde{\varphi} : A \to \Bbb L(l_2(B \otimes \mathcal K))$
is an equivariant $*$-homomorphism modulo $\Bbb K(l_2(B \otimes
\mathcal K))$. By identifying $\Bbb L(l_2(B\otimes \mathcal K))$
with $M(B\otimes \mathcal K)$, $\Bbb K(l_2(B\otimes \mathcal K))$
with $B\otimes \mathcal K$ and the quotient $ \Bbb L(l_2(B\otimes
\mathcal K))/\Bbb K(l_2(B\otimes \mathcal K))$ with $Q(B\otimes
\mathcal K)$, we can consider $\widetilde{\varphi}$ as a map
$\widetilde{\varphi} : A \to M(B\otimes \mathcal K)$ with the property
that $q_{B \otimes \mathcal K} \circ \widetilde{\varphi} \in
\Hom_G(A,Q(B\otimes \mathcal K))$.

\begin{lemma}\label{LLL2} The class of $q_{B \otimes \mathcal K} \circ
\widetilde{\varphi}$ in $\Ext(A,B \otimes \mathcal K)_h$ is independent
of the choice of unit sequence, subject to the conditions 0)-8), and of
the chosen discretization, subject to condition 9), and depends only on
the class $[\varphi]$ of $\varphi$ in $[[A,Q(B)_G \otimes \mathcal K]]$.
\begin{proof} Let $\{v_n\}$ be another unit sequence satisfying 0)-8).
There is then a unit sequence $\{w_n\}$ in $B \otimes \mathcal K$ such
that $w_nv_n = v_n, w_n u_n = u_n$ for all $n$. Connect $u_0$ to $w_0$
by a straight line, then $u_1$ to $w_1$ by a straight line, etc. This
gives a path $\{w_n^t\}_{t \in [0,1[}$ of unit sequences. For each
$t \in [0,1[$ we get then a map $\mu_t : A \to M(B \otimes \mathcal K)$
such that $q_{B \otimes \mathcal K} \circ \mu_t \in \Hom_G(A,Q(B \otimes
\mathcal K))$ and $[q_{B\otimes \mathcal K} \circ \mu_0] =
[q_{B\otimes \mathcal K} \circ \widetilde{\varphi}]$ in $\Ext(A,B
\otimes \mathcal K)$. Let $\delta : A \to M(B \otimes \mathcal K)$
be the map obtained from $\varphi$ as $\widetilde{\varphi}$ was, but
by using $\{w_n\}$ instead of $\{u_n\}$. Then $\lim_{t \to 1} \mu_t(a)
= \delta(a)$ in the strict topology for all $a \in A$, and
$$
\lim_{t \to 1} \mu_t(a)\mu_t(b) - \mu_t(ab) = \delta(a)\delta(b) -
\delta(ab),
$$
$$
\lim_{t \to 1} \mu_t(a + \lambda b) - \mu_t(a) - \lambda\mu_t(b) =
\delta(a +b) - \delta(a) - \lambda \delta(b) ,
$$
$$
\lim_{t \to 1} \mu_t(a^*) - \mu_t(a)^* = \delta(a^*) - \delta(a)^* ,
$$
$$
\lim_{t \to 1} \mu_t( g \cdot a) - g \cdot \mu_t(a) = \delta(g \cdot a) -
g \cdot \delta(a) ,
$$
in norm for all $a,b \in A, \lambda \in \Bbb C, g \in G$. Hence
$[q_{B \otimes \mathcal K} \circ \delta] = [q_{B \otimes \mathcal K}
\circ \widetilde{\varphi}]$ in $\Ext(A,B \otimes \mathcal K)_h$. The
same argument with the unit sequence $\{u_n\}$ replaced by $\{v_n\}$
shows that the class of $[q_{B \otimes \mathcal K} \circ
\widetilde{\varphi}]$ in $\Ext(A,B \otimes \mathcal K)_h$ is
independent of the choice of unit sequence. Once this is established
it is clear that a homotopy of asymptotic homomorphisms $A \to Q(B)_G
\otimes \mathcal K$ gives rise, by an appropriate choice of unit sequence,
to a homotopy which shows that $[q_{B \otimes \mathcal K} \circ
\widetilde{\varphi}] \in \Ext(A,B \otimes \mathcal K)_h$ only depends
on the homotopy class of $\varphi$. That $[q_{B \otimes \mathcal K}
\circ \widetilde{\varphi}]$ is also independent of the discretization
and only depends on the homotopy class of $\varphi$ follows in the same
way as in Lemma 5.3 and Lemma 5.4 of \cite{MT1}.

\end{proof}
\end{lemma}

It follows that we have the desired map $\alpha : [[A,Q(B)_G \otimes
\mathcal K]] \to \Ext(A,B \otimes \mathcal K)_h$ which is easily seen
to be a semi-group homomorphism.

\begin{lemma}\label{LLL3} Let $\varphi : SA \to Q(B) \otimes \mathcal K$
be an equivariant $*$-homomorphism which we consider as a (constant)
asymptotic homomorphism. Let $X$ be a compact subset with dense span
in $SA$ and choose a unit sequence $\mathcal U = \{u_n\}$ in $B \otimes
\mathcal K$ such that
\eq{\label{eee1}
\|\sqrt{u_n - u_{n-1}} \chi(\varphi(a)) - \chi(\varphi(a))\sqrt{u_n -
u_{n-1}} \| < 2^{-n} }
for all $a \in X$ and
\eq{\label{eee2}
\sum_{j=1}^{\infty} \| g \cdot \sqrt{u_j - u_{j-1}} -
\sqrt{u_j - u_{j-1}} \|^2 < \infty }
for all $g \in G$. Then $[q_{B \otimes \mathcal K} \circ
\widetilde{\varphi}] = [\iota \circ \varphi]$ in
$\Ext(SA,B \otimes \mathcal K)$, where $\iota : Q(B)_G \otimes \mathcal
K \to Q(B\otimes \mathcal K)_G$ is the natural embedding.
\begin{proof} $\widetilde{\varphi}$ has the form $\widetilde{\varphi}(a)
= P_{\mathcal U} ( \sum_{j = 0}^{\infty} \chi(\varphi(a))e_{jj})
P_{\mathcal U}$. Let $V \in \Bbb L(l_2(B \otimes \mathcal K))$ be
the partial isometry defining $P_{\mathcal U}$ and note that
$g \cdot V - V \in \Bbb K(l_2(B \otimes \mathcal K))$ for all
$g \in G$ because of (\ref{eee2}). Thus
$$
\left ( \begin{matrix} V & 1 - VV^* \\ 1 - V^*V & - V^*
\end{matrix} \right )
$$
is a unitary in $M_2(\Bbb L(l_2(B \otimes \mathcal K)))$ which is
$G$-invariant modulo $M_2(\Bbb K(l_2(B \otimes \mathcal K)))$ and
satisfies that
$$
\left ( \begin{matrix} V & 1 - VV^* \\ 1 - V^*V & - V^* \end{matrix}
\right )  \left( \begin{matrix} \widetilde{\varphi}  & 0 \\ 0 & 0
\end{matrix} \right ) \left ( \begin{matrix} V^* & 1 - V^*V \\ 1 -
VV^* & - V \end{matrix} \right ) = \left ( \begin{matrix} \varphi_0
& 0 \\ 0 & 0 \end{matrix} \right ) ,
$$
where $\varphi_0(a) = (\sqrt{u_0} \chi(\varphi(a))\sqrt{u_0} + \sum_{j
= 1}^{\infty} \sqrt{u_j - u_{j-1}}\chi(\varphi(a))\sqrt{u_j - u_{j-1}})
e_{00}$. Thanks to (\ref{eee1}) we have that
$$
\sum_{j = 1}^{\infty} \|\sqrt{u_j - u_{j-1}}\chi(\varphi(a))\sqrt{u_j -
u_{j-1}} - (u_j - u_{j-1})\chi(\varphi(a))\| < \infty
$$
for all $a \in X$. Since $\sum_{j = 1}^{\infty} (u_j - u_{j-1})
\chi(\varphi(a)) + u_0 \chi(\varphi(a)) = \chi(\varphi(a))$
(with convergence in the strict topology) we find that $\varphi_0(a) =
\chi(\varphi(a))e_{00}$ modulo $\Bbb K(l_2(B \otimes \mathcal K))$
for all $a \in X$, and hence in fact for all $a \in SA$. This proves
the lemma.

\end{proof}
\end{lemma}

\section{The main results}

Since $A$ is separable, $[[SA, X \otimes \mathcal K]] = \varinjlim_D
[[SA,D \otimes \mathcal K]]$ for any $G$-algebra $X$, when we take the
limit over all separable $G$-subalgebras $D$ of $X$. It follows from
\cite{DL} that the suspension map $S : [[SA, X \otimes \mathcal K]]
\to [[S^2A, SX \otimes \mathcal K]]$ is an isomorphism.\footnote{Dadarlat
and Loring did not consider the equivariant theory in \cite{DL}, but it
is easy to check that their arguments carry over unchanged.} Hence
$[[SA, - \otimes \mathcal K]]$ is a homotopy invariant and half-exact
functor on the category of $G$-algebras (and not only separable
$G$-algebras). There is therefore a map
$$
\partial : [[SA, SQ(B)_G\otimes \mathcal K]] \to [[SA,B \otimes \mathcal K]]
$$
arising as the boundary map coming from the extension (\ref{ex1}), cf.
e.g. \cite{GHT}. Well-known arguments from the K-theory of $C^*$-algebras,
cf. \cite{Bl}, show that $[[SA, SM(B)_G \otimes \mathcal K]] =
[[SA, M(B)_G \otimes \mathcal K]] = 0$, so the six-terms exact sequence
obtained by applying $[[SA, - \otimes \mathcal K]]$ to (\ref{ex1}) shows
that $\partial$ is an isomorphism. For any $G$-algebra $D$ we let $ s :
D \to D \otimes \mathcal K$ be the stabilising $*$-homomorphism given by
$s(d) = d \otimes e$ for some minimal projection $e \in \mathcal K$.
Since $B$ is weakly stable there is an equivariant $*$-isomorphism
$\gamma_0 : B \otimes \mathcal K \to B$ such that $s \circ \gamma_0 :
B \otimes \mathcal K \to B \otimes \mathcal K$ is equivariantly
homotopic to $\id_{B \otimes \mathcal K}$. Let $\gamma : Q(B \otimes
\mathcal K)_G \to Q(B)_G$ the $*$-isomorphism induced by $\gamma_0$.

\begin{lemma}\label{LC1} The composition of the maps
\eQ{\begin{xymatrix}{
[[S^2A, B \otimes \mathcal K]] \ar[r]^-{\partial^{-1}} & [[S^2A ,
SQ(B)_G \otimes \mathcal K]] \\
\ar[r]^-{S^{-1}} &  [[SA,Q(B)_G \otimes \mathcal K]] \ar[r]^-{\alpha}
& \Ext(SA, B \otimes \mathcal K)_h  \ar[r]^-{CH} &  [[S^2A,B \otimes
\mathcal K]] } \end{xymatrix} }
is the identity.
\begin{proof} We are going to use Theorem 2.3 of \cite{H-LT}.\footnote{The
equivariant theory was not explicitly considered in \cite{H-LT}, but all
arguments carry over unchanged.} Let $ x = s_*([\id_{SB}]) \in [[SB,SB
\otimes \mathcal K]]$, where $[\id_{SB}] \in [[SB,SB]]$ is the element
represented by the identity map of $SB$ and $s : SB \to SB \otimes
\mathcal K$ is the stabilising $*$-homomorphism. By Theorem 2.3 of
\cite{H-LT} it suffices to identify the image of $x$ under the
Bott-periodicity isomorphism $[[SB, SB\otimes \mathcal K]] \simeq
[[S^2B, B \otimes \mathcal K]]$  and show that the image of that
element is not changed under the map we are trying to prove is always
the identity. This is what we do. Under the isomorphism $[[SB,SB \otimes
\mathcal K]] \simeq [[S^2B, B \otimes \mathcal K]] $, coming from
Bott-periodicity, the image of $x$ is represented by the asymptotic
homomorphism $S^2B \to B \otimes \mathcal K$ arising by applying the
Connes-Higson construction to the Toeplitz extension tensored with $B$ :
\eq{\label{ex2}\begin{xymatrix}{
0 \ar[r]   &  B \otimes \mathcal K  \ar[r] &  T_0 \otimes B  \ar[r]  &
SB \ar[r]  & 0 .}
\end{xymatrix}}
In other words, if $\varphi : SB \to Q(B \otimes \mathcal K)$ is the
Busby invariant of (\ref{ex2}) the image of $x$ in $[[S^2B,B \otimes
\mathcal K]]$ is $[CH(\varphi)]$. For each separable $G$-subalgebra $D
\subseteq Q(B)_G$ we let $\iota_D : D \to Q(B)_G$ denote the inclusion.
Then the boundary map $\partial : [[S^2B , SQ(B)_G \otimes \mathcal K]]
\to [[S^2B, B \otimes \mathcal K]]$ is given by
$$
\partial (z) \ = \ \lim_D [CH(\iota_D) \otimes \id_{\mathcal K}] \bullet z ,
$$
where $\bullet$ denote the composition product in $E$-theory. Hence
$\partial^{-1}[CH(\varphi)]$ is the element $z \in [[S^2B, SQ(B)_G
\otimes \mathcal K]]$ with the property that
$$
\lim_D[CH(\iota_D) \otimes \id_{\mathcal K}] \bullet z = [CH(\varphi)]
$$
for all large enough $D$. Let $\iota : Q(B)_G \otimes \mathcal K \to
Q(B \otimes \mathcal K)_G$ be the natural embedding. By the naturality
of the Connes-Higson construction,
$$
[CH(\iota_D) \otimes \id_{\mathcal K}]  \bullet S([s \circ \gamma
\circ \varphi]) = [CH(\iota \circ s \circ \gamma \circ \varphi)]
$$
for all separable $G$-subalgebras $D \subseteq Q(B)_G$ which contains
$\gamma \circ \varphi(SB)$. Since $s \circ \gamma_0$ is equivariantly
homotopic to the identity map, we have that
$$
 [CH(\iota \circ s \circ \gamma \circ \varphi)] = (s \circ \gamma_0)_*
[CH(\varphi)] = [CH(\varphi)],
$$
so we conclude that $\partial^{-1}[CH(\varphi)] = S([s \circ \gamma
\circ \varphi])$. Hence $\alpha \circ S^{-1} \circ \partial^{-1}
[CH(\varphi)] = [ \iota \circ s \circ \gamma \circ \varphi]$ by
Lemma \ref{LLL3}. Thus the image of $[CH(\varphi)]$ in $[[S^2B, B
\otimes \mathcal K]]$ under the composite map is $CH[\iota \circ s
\circ \gamma \circ \varphi] = (s \circ \gamma_0)_*[CH(\varphi)] =
[CH(\varphi)]$. The proof is complete.
\end{proof}
\end{lemma}

\begin{lemma}\label{CC1} Let $\lambda \in \Ext(SA,B \otimes \mathcal K)$.
Then $\varphi = s \circ \gamma \circ \lambda$ is an equivariant
$*$-homomorphism $\varphi : SA \to Q(B)_G\otimes \mathcal K$ such that
$\alpha[\varphi] = s_* \circ \gamma_* [\lambda]$ in $\Ext(SA, B \otimes
\mathcal K)_h$ and such that $[\varphi] = 0$ in $[[SA,Q(B)_G\otimes
\mathcal K]]$ implies that $[\lambda] = 0$ in $\Ext(SA,B \otimes
\mathcal K)$.
\begin{proof} If $[\varphi] = 0$ in $[[SA, Q(B)_G\otimes \mathcal K]]$,
there is a path $\mu^t, t \in [0,1]$, of asymptotic homomorphisms
$SA \to Q(B)_G \otimes \mathcal K$ such that $\mu^0 = \varphi$ and
$\mu^1 = 0$ and a unit sequence $\mathcal U = \{u_n\}$ in $B \otimes
\mathcal K$ such that
\eq{\label{e6}
q_{B \otimes \mathcal K} \circ \widetilde{\mu^t} , \ t \in [0,1], }
connects $q_{B \otimes \mathcal K}\circ \widetilde{\varphi}$ to $0$.
By Theorem \ref{TT1} we may assume that $\mu$ is an equi-homotopy and
it is then easy to see that (\ref{e6}) is a strong homotopy. By Lemma
\ref{LL1} we conclude from this that $[q_{B \otimes \mathcal K}
\circ \widetilde{\varphi}] = 0$ in $\Ext(SA, B \otimes \mathcal K)$.
But $[q_{B \otimes \mathcal K} \circ \widetilde{\varphi}] = [\varphi]$
in $\Ext(SA,B \otimes \mathcal K)$ by Lemma \ref{LLL3}.
Hence $\alpha[\varphi] = s_* \circ \gamma_*[\lambda]$ in
$\Ext(SA, B\otimes \mathcal K)_h$ and $[\varphi] = 0 \Rightarrow
s_* \circ  \gamma_*[\lambda] = 0$ in $\Ext(SA, B \otimes \mathcal K)$.
To complete the proof it suffices to show that $s_* \circ \gamma_* :
\Ext(SA,B \otimes \mathcal K) \to \Ext(SA, B \otimes \mathcal K)$ is
injective. However, $\gamma$ is an equivariant $*$-isomorphism and
therefore $\gamma_*$ is an isomorphism. The injectivity of $s_* :
\Ext(SA,B) \to \Ext(SA,B \otimes \mathcal K)$ follows from
the weak stability of $B$ : There is a $G$-invariant isometry
$V \in M(B \otimes \mathcal K)$ such that $x \mapsto V^* s(x)V$
is an equivariant $*$-automorphism $B\otimes \mathcal K \to B \otimes
\mathcal K$ and $s(x) = \Ad V(V^*s(x) V)$. Since $\Ad V$ induces the
identity map on $\Ext(SA,B \otimes \mathcal K)$ we see that $s_* :
\Ext(SA, B) \to \Ext(SA, B \otimes \mathcal K)$ is an isomorphism.

\end{proof}
\end{lemma}

\begin{lemma}\label{CH1-1} The map $CH : \Ext(SA,B) \to [[S^2A,B]]$ is
injective.
\begin{proof}
Consider an extension $\lambda \in \Ext(SA, B \otimes \mathcal K)$ and
assume that $[CH(\lambda)] = 0$ in $[[S^2A,B \otimes \mathcal K]]$. With
the notation from Lemma \ref{CC1} we find that $CH \circ \alpha [\varphi]
= CH[s \circ \gamma \circ \lambda]  = s_* \circ \gamma_*[CH(\lambda)] = 0$.
But then Lemma \ref{LC1} implies that $[\varphi] = 0$ in $[[SA, Q(B)_G
\otimes \mathcal K]]$. By Lemma \ref{CC1} this yields the conclusion
that $[\lambda] = 0$ in $\Ext(SA, B \otimes \mathcal K)$. Thus $CH :
\Ext(SA,B \otimes \mathcal K) \to [[S^2A,B \otimes \mathcal K]]$ is
injective. But $B$ is weakly stable so the result follows.
\end{proof}
\end{lemma}

\begin{cor}\label{coinc} $\Ext(SA,B) = \Ext(SA,B)_h$.
\end{cor}

\bigskip

The surjectivity of $CH : \Ext(SA,B) \to [[S^2A,B]]$ follows from
Lemma \ref{LC1}. Furthermore, it follows from Lemma \ref{CH1-1}
that $\alpha$ is well-defined as a map $\alpha : [[SA, Q(B)_G \otimes
\mathcal K]] \to \Ext(SA, B \otimes \mathcal K)$ and then Lemma \ref{LC1}
tells us that
$$
CH^{-1} = \alpha \circ S^{-1} \circ \partial^{-1}.
$$
Another description of $CH^{-1}$ can be obtained from \cite{MT2}.
The crucial construction for this is the map $E$ which was considered
in \cite{MT1} and \cite{MT2}, inspired by \cite{MM} and \cite{MN}.
However, in \cite{MT1} and \cite{MT2} we only defined $E$ as a map
into homotopy classes of extensions, so to see that the $E$-construction
can also invert the CH-map of Lemma \ref{CH1-1} we must show that it is
well-defined as a map from homotopy classes of asymptotic homomorphisms
to stable unitary equivalence classes of extensions. Let us therefore
review the construction.

Given an equicontinuous asymptotic homomorphism
$\varphi = \{\varphi_t\}_{t \in [1,\infty)} : A \to B$ we choose a
discretization $\{\varphi_{t_i}\}_{i \in \Bbb N}$ such that
$\lim_{i \to \infty} t_i = \infty $ and $\lim_{i \to \infty}
\sup_{t  \in [t_i,t_{i+1}]} \|\varphi_t(a) - \varphi_{t_i}(a)\| = 0$
for all $a \in A$. Since $G$ is $\sigma$-compact (and $\varphi$
equicontinuous) we can also arrange that
$$
\lim_{i \to \infty} \sup_{t  \in [t_i,t_{i+1}]} \sup_{g \in K} \|
g \cdot \varphi_t(a) - \varphi_{t}(g \cdot a)\| = 0
$$
for all $a \in A$ and all compact subsets $K \subseteq G$. To define
from such a discretization a map $\bold\Phi : A \to \Bbb L(l_2(\Bbb Z)
\otimes B)$ we introduce the standard matrix units $e_{ij},
i,j \in \Bbb Z$, which act on the Hilbert $B$-module $l_2(\Bbb Z)\otimes B$
in the obvious way. Then
$$
\bold\Phi (a) = \sum_{i \geq 1} \varphi_{t_i}(a) e_{ii}
$$
defines a map $\bold\Phi : A \to \Bbb L(l_2(\Bbb Z) \otimes B)$.
As in the proof of Lemma \ref{crux} we can define a representation of
$G$ on $l_2(\Bbb Z) \otimes B$ and in this way obtain a representation
of $G$ as automorphisms of $\Bbb L(l_2(\Bbb Z) \otimes B)$. Since $B$
is weakly stable we can identify $B$ with $\Bbb K(l_2(\Bbb Z) \otimes B))$,
the $B$-compact operators in $\Bbb L(l_2(\Bbb Z) \otimes B)$. Observe
that $\bold\Phi$ is then an equivariant $*$-homomorphism modulo $B$.
Furthermore, $\bold\Phi(a)$ commutes modulo $B$ with the two-sided shift
$T = \sum_{j \in \Bbb Z} e_{j,j+1}$ which is $G$-invariant. So we get in
this way a $G$-extension $E(\varphi) : A \to Q(B) =  \Bbb L(l_2(\Bbb Z)
\otimes B)/\Bbb K(l_2(\Bbb Z) \otimes B)$ such that
$$
E(\varphi)(f \otimes a) = f(\underline{T})\underline{\bold\Phi(a)}
$$
for all $f \in C(\Bbb T), a \in A$. Here and in the following we denote
by $\underline{S}$ the image in $Q(B) =
\Bbb L(l_2(\Bbb Z) \otimes B)/\Bbb K(l_2(\Bbb Z) \otimes B)$
of an element $S \in \Bbb L(l_2(\Bbb Z) \otimes B)$.

\begin{lemma}\label{LL2}
$E( \varphi)$ is a semi-invertible $G$-extension, and up to stable
unitary equivalence it does not depend on the chosen discretization of
$\varphi$.
\begin{proof} Consider another discretization $(\varphi_{s_i})_{i \in
\Bbb N}$ of $\varphi$ and define $\bold\Psi : A \to \Bbb L(l_2(\Bbb Z)
\otimes B)$ by
$$
\bold\Psi(a) = \sum_{i \leq 0} \varphi_{s_{-i +1}}(a) e_{ii} .
$$
There is then a $G$-extension $-E(\varphi) : C(\Bbb T) \otimes A \to
\Bbb L(l_2(\Bbb Z) \otimes B)/\Bbb K(l_2(\Bbb Z) \otimes B)$ such
that $-E(\varphi)(f \otimes a) = f(\underline{T})\underline{\bold\Psi(a)}$.
It suffices to show that $-E(\varphi) \oplus E(\varphi)$ is unitarily
equivalent to an asymptotically split $G$-extension. Define
$\Lambda : A \to \Bbb L(l_2(\Bbb Z) \otimes B)$ such that
$$
\Lambda(a) = \sum_{i \geq 1} \varphi_{t_i}(a)e_{ii} + \sum_{i \leq 0}
\varphi_{s_{- i+1}}(a)e_{ii} .
$$
There is then a $G$-extension $\pi_0 : C(\Bbb T) \otimes A \to
\Bbb L(l_2(\Bbb Z) \otimes B)/\Bbb K(l_2(\Bbb Z) \otimes B)$ such
that $\pi_0( f \otimes a) = f(\underline{T})\underline{\Lambda(a)}$.
$-E(\varphi) \oplus E(\varphi)$ is clearly unitarily equivalent (via a
$G$-invariant unitary) to $\pi_0 \oplus 0$, so it suffices to show
that $\pi_0$ is asymptotically split. For each $n$ we define
$\Lambda_n : A \to \Bbb L(l_2(\Bbb Z) \otimes B)$ by
\eQ{\spl{
&\Lambda_n(a) = \\
&\sum_{i > n} \varphi_{t_i}(a)e_{ii} + \sum_{1 \leq i \leq n}
\varphi_{t_n}(a) e_{ii} + \sum_{\{i \leq 0 : \ s_{-i+1} \leq t_n\}}
\varphi_{t_n}(a) e_{ii} +  \sum_{\{i \leq 0 : \ s_{-i+1} > t_n\}}
\varphi_{s_i}(a) e_{ii}. }}
Then $\{\Lambda_n\}_{n \in \Bbb N}$ is a discrete asymptotic
homomorphism such that $\lim_{n \to \infty} \|\Lambda_n(a) -
\Lambda_{n+1}(a)\| = 0$, $\lim_{n \to \infty} \|g \cdot \Lambda_n(a) -
\Lambda_n(g \cdot a)\| = 0, g \in G$, $\lim_{n \to \infty} \|T
\Lambda_n(a) - \Lambda_n(a)T\| = 0$ and $\Lambda_n(a) = \Lambda(a)$
modulo $\Bbb K(l_2(\Bbb Z)\otimes B)$. By convex interpolation and
an obvious application of the $C^*$-algebra
$$
\{ f \in C_b([1,\infty), M(B))  \ : \ q_B(f(t)) = q_B(f(1)),
t \in [1,\infty) \}/C_0([1,\infty),B)
$$
we get an asymptotic homomorphism $(\pi_t)_{t \in [1,\infty)} :
C(\Bbb T) \otimes A \to M(B) = \Bbb L(l_2(\Bbb Z) \otimes B)$ such
that $\pi_0 = q_B \circ \pi_t$ for all $t$.
\end{proof}
\end{lemma}

Theorem \ref{TT1} and Lemma \ref{LL2} in combination show that there
is group homomorphism $E : [[SA,B]] \to \Ext(C(\Bbb T) \otimes SA,B)$
such that $E[\varphi] = [E(\varphi)]$ for any equicontinuous asymptotic
homomorphism $ \varphi : SA \to B$. By pulling extensions back along the
inclusion $S^2A  \subseteq C(\Bbb T) \otimes SA$ we can also consider
$E$ as a map $E : [[SA,B]] \to \Ext(S^2A,B)$. Let $\chi : SA \to S^3M_2(A)$
be a $*$-homomorphism which is invertible in KK-theory. By weak stability
of $B$ there is also an isomorphism $\beta : [[S^2A,B]] \to
[[S^2M_2(A),B]]$. Let $\xi : S^2 \to \mathcal K$ be the asymptotic
homomorphism which arises from the Connes-Higson construction applied
to the Toeplitz extension. By changing $\chi$ 'by a sign' we may assume
that the composite map
\eQ{\begin{xymatrix}{
[[S^2A, B]] \ar[r]^-{\beta}  & {[[S^2M_2(A),B]] \ \ \ \ar[r]^-{{[\varphi]
\mapsto [\xi \otimes \varphi]}} \  } &   \ \ \ [[S^4M_2(A),B]]
\ar[r]^-{(S\chi)^*} & [[S^2A,B]]  }
\end{xymatrix} }
is the identity. Consider the diagram
\eQ{\begin{xymatrix}{
\Ext(SA,B) \ar[d]_-{CH}  & \Ext(S^3M_2(A),B)  \ar[l]_-{\chi^*}
\ar[d]^-{CH} \\
[[S^2A,B]]   \ar[ur]^-{E \circ \beta} &  [[S^4M_2(A),B]]
\ar[l]^-{(S\chi)^*} . }
\end{xymatrix} }
The square commutes by the naturality of the Connes-Higson construction,
and it follows from Lemma 2.3 of \cite{MT2} (or Lemma 5.5 of \cite{MT1})
that $(S \chi)^* \circ CH \circ E \circ \beta = \id$. We conclude
therefore that $CH \circ \chi^* \circ E \circ \beta = \id$. We have
now obtained our main results :

\bigskip

\begin{thm}\label{MAIN1}
Let $A$ and $B$ be separable $G$-algebras, $B$ weakly stable.
$CH : \Ext(SA,B) \to [[S^2A,B]]$ is an isomorphism with inverse
$\chi^* \circ E \circ \beta$.
\end{thm}

\bigskip

It follows, of course, that the bifunctor $\Ext(SA,B)$ has the same
properties as $E$-theory, such as excision and Bott periodicity in
both variables, for example.

\bigskip

\begin{thm}\label{MAIN}
Let $A$ and $B$ be separable $G$-algebras, $B$ weakly stable,
and let $\varphi, \psi : SA \to Q(B)$ be two $G$-extensions.
The following conditions are equivalent :
\begin{enumerate}
\item[1)] $[\varphi] = [\psi]$ in $\Ext(SA,B)$ (i.e. $\varphi$ and
$\psi$ are stably unitarily equivalent).
\item[2)] $\varphi \oplus 0$ and $\psi \oplus 0$ are strongly homotopic.
\item[3)] $\varphi$ and $\psi$ are homotopic.
\end{enumerate}
\begin{proof} 1) $\Rightarrow$ 2): Assuming 1) there is an
asymptotically split extension $\lambda$ such that $\varphi \oplus
\lambda$ and $\psi \oplus \lambda$ are unitarily equivalent. By
Lemma 6.1 of \cite{Th1} this implies that $\varphi \oplus \lambda
\oplus 0$ and $\psi \oplus \lambda \oplus 0$ are strongly homotopic.
Then $\varphi \oplus \lambda \oplus (\lambda \circ \alpha) \oplus 0$
and $\psi \oplus \lambda \oplus (\lambda \circ \alpha) \oplus 0$ are
also strongly homotopic, where $\alpha \in \Aut SA$ inverts the
orientation of the suspension. 2) follows by observing that
$\lambda \oplus (\lambda \circ \alpha)$ is strongly homotopic to $0$.
2) $\Rightarrow$ 3) follows because an invariant isometry in
$M(B)$ can be connected to $1$ via a strictly continuous path of
$G$-invariant isometries, cf. e.g. Lemma 3.3 2) of \cite{Th1}.
3) $\Rightarrow$ 1) follows from Lemma \ref{CH1-1}.
\end{proof}
\end{thm}

\bigskip

\begin{remark} It is easy to extend Theorem \ref{MAIN1} and
Theorem \ref{MAIN} to the case where $B$ is only $\sigma$-unital
(i.e. contains a strictly positive element). In fact, it suffices
to observe that
$$
\Ext(SA,B) \simeq \varinjlim_D \Ext (SA, D),
$$
where we take the limit over all weakly stable separable
$G$-subalgebras $D$ of $B$ with the property that $D$ contains a
positive element which is strictly positive in $B$.
\end{remark}

\section{$K$-homology}

It follows from Theorem \ref{MAIN1} and Theorem \ref{MAIN} that
$\Ext(SA,B) = [[S^2A,B]]$ can also be identified with the homotopy
classes of equivariant $*$-homomorphisms $\psi : SA \to Q(B)$ with
the property that $\psi \oplus 0$ is strongly homotopic to $\psi$.
As a consequence we conclude that
$$
[[S^2A,B]] \ \simeq \ \varinjlim_n [SA, Q(B)\otimes M_n(\Bbb C)],
$$
where $[\cdot \ , \ \cdot]$ denotes homotopy classes of equivariant
$*$-homomorphisms. In the important special case where $B = \mathcal K$,
and the group $G$ is trivial, we can even do better. Let $Q$ denote the
Calkin algebra, $Q = \Bbb L(l_2)/\Bbb K(l_2)$.

\begin{lemma}\label{lllem} Let $\varphi : SA \to Q$ be a $*$-homomorphism.
There is then an isometry $V \in \Bbb L(l_2)$ with infinite dimensional
co-kernel and a $*$-homomorphism $\varphi_0 : SA \to Q$ such that
$\varphi$ is homotopic to $\Ad q_{\mathcal K}(V) \circ \varphi_0$.
\begin{proof} We may assume that $\varphi$ is not homotopic to $0$.
Let $\iota_i : SA \to SA, i = 1,2$, be $*$-homomorphisms with orthogonal
ranges, both homotopic to the identity map. Then $\varphi \circ \iota_1$
and $\varphi \circ \iota_2$ are homotopic to $\varphi$, and in
particular non-zero. Let $a$ be a non-zero positive element in the
range of $\varphi \circ \iota_2$ and let $b \in \Bbb L(l_2)$ be a
positive lift of $a$. By spectral theory $b\Bbb L(l_2)b$ contains a
projection $E$ with non-zero image in $Q$. Since
$(1 - q_{\mathcal K}(E))x = x$ for all $x \in \varphi \circ \iota_1(SA)$,
we conclude that $1 - q_{\mathcal K}(E)$ is non-zero in $Q$.
It follows that there is an isometry $V$ with infinite dimensional
co-kernel such that $VV^* = 1 -E$. Set $\varphi_0 =
\Ad q_{\mathcal K}(V^*) \circ \varphi \circ \iota_1$.

\end{proof}
\end{lemma}

\begin{thm}\label{K-homology}
Let $A$ be a separable $C^*$-algebra. Then $E(A, \Bbb C)$ is naturally
isomorphic to the group $[SA,Q]$ of homotopy classes of $*$-homomorphisms
from $SA$ to $Q$.
\begin{proof}  It follows from Lemma \ref{lllem} that $\varphi$ is
strongly homotopic to $\varphi \oplus 0$ for any $*$-homomorphism
$\varphi : SA \to Q$. By using that the unitary group of $\Bbb L(l_2)$
is norm-connected, it follows from this and Theorem \ref{MAIN} that
$\Ext(SA,\mathcal K)$ is naturally isomorphic to $[SA,Q]$. Since
$\Ext(SA,\mathcal K)$ is naturally isomorphic to $E(A,\Bbb C)$ by
Theorem \ref{MAIN1}, this completes the proof.
\end{proof}
\end{thm}

A weak version of Theorem \ref{K-homology} was conjectured by
Rosenberg in \cite{R}.

\vspace{2cm}
\parbox{7cm}{V. M. Manuilov\\
Dept. of Mech. and Math.,\\
Moscow State University,\\
Moscow, 119899, Russia\\
e-mail: manuilov@mech.math.msu.su
}
\hfill
\parbox{6cm}{K. Thomsen\\
Institut for matematiske fag,\\
Ny Munkegade, 8000 Aarhus C,\\
Denmark\\
e-mail: matkt@imf.au.dk
}

\end{document}